\documentclass[titlepage]{amsart}
\usepackage[a4paper,top=2cm,bottom=2cm,left=2cm,right=2cm]{geometry}
\usepackage{lipsum}
\makeatletter
\g@addto@macro{\endabstract}{\@setabstract}
\newcommand{\authorfootnotes}{\renewcommand\thefootnote{\@fnsymbol\c@footnote}}%
\makeatother
\usepackage{microtype} 
\usepackage[osf,sc]{mathpazo} 
\usepackage{amsmath,amssymb,amsthm}    
\usepackage{amsfonts}
\usepackage{mathrsfs}
\usepackage{mathtools}
\usepackage{stmaryrd}  
\usepackage{arydshln}
\usepackage[autostyle,italian=guillemets]{csquotes}
\usepackage[backend=biber,style=numeric-comp]{biblatex}

\providecommand{\gtrless}{
	\mathrel{
		\smash{
			\vcenter{
				\offinterlineskip 
				\ialign{
					\hfil##\hfil\cr 
					$>$\cr 
					\noalign{\kern-.3ex}
					$<$\cr 
				}
			}
		}
		\vphantom{>}
	}
}

\usepackage[figuresright]{rotating}

\begin{document}
	\begin{center}
		\LARGE
		A stochastic approach to delays optimization for narrowband transmit beam pattern in medical ultrasound\par \bigskip
		\normalsize
		\authorfootnotes
		Chiara Razzetta\footnote{razzetta@dima.unige.it - Corresponding Author}\textsuperscript{1}, Valentina Candiani \footnote{candiani@dima.unige.it - Corresponding Author}\textsuperscript{1}, Marco Crocco\footnote{marco.crocco@esaote.com}\textsuperscript{2}, Federico Benvenuto\footnote{benvenuto@dima.unige.it}\textsuperscript{1}
		\par \bigskip
		\textsuperscript{1}DIMA - Università degli Studi di Genova, Via Dodecaneso 35, I-16146 Genoa, Italy\par 
		\textsuperscript{2}Esaote S.P.A. Via E. Melen 77, I-16152 Genoa, Italy
	\end{center}

\begin{abstract}
Ultrasound imaging is extensively employed in clinical settings due to its non-ionizing nature and real-time capabilities. The beamformer represents a crucial component of an ultrasound machine, playing a significant role in shaping the ultimate quality of the reconstructed image.
Therefore, Transmit Beam Pattern (TBP) optimization is an important task in medical ultrasound, but state-of-the-art TBP optimization has well-known drawbacks like non-uniform beam width over depth, presence of significant side lobes, and quick energy drop out after the focal depth. 
To overcome these limitations, we developed a novel optimization approach for TBP by focusing the analysis on its narrowband approximation, particularly suited for Acoustic Radiation Force Impulse (ARFI) elastography, and considering transmit delays as free variables instead of linked to a specific focal depth. We formulate the problem as a non linear Least Squares problem to minimize the difference between the TBP corresponding to a set of delays and the desired one, modeled as a 2D rectangular shape elongated in the direction of the beam axis. In order to quantitatively evaluate the results, we define three quality metrics based on main lobe width, side lobe level, and central line power.
Results obtained by our synthetic software simulation show that the main lobe width is considerably more intense and uniform over the whole depth range with respect to classical focalized Beam Patterns, and our optimized delay profile results in a combination of standard delay profiles at different focal depths.
The application of the proposed method to ARFI elastography shows improvements in the concentration of the ultrasound energy along a desired axis.
\end{abstract}

\vspace{2pc}

\section{Introduction}
Medical Ultrasound Imaging is the most widespread real-time non-invasive imaging system: it exploits the ability of human tissue to reflect ultrasound signals.
In particular, different ultrasound modalities process the reflected echoes, generating either morphological images or measures of tissue properties. 
The most widely recognized biomedical ultrasound technique, known as B-mode, produces the standard black and white images by means of several transmission and reception phases, resulting in a map of acoustic tissue impedance.
Over the last decades, advancements in machine capabilities, image quality, and computational power have paved the way to innovative strategies, techniques, and methods.
In particular, several quantitative modalities have been developed, allowing to measure physical parameters of tissues with diagnostic significance (Doppler modalities~\cite{Brody1974} and elastography~\cite{Doherty2013} can be cited among them). 
The performance of each modality significantly depends on how the transmission and reception phases are tuned, reflecting both on the image quality and the accuracy of measurements. 
The two phases are characterized by different parameters such as pulse shape, central frequency, transmission focus, transmit and receive apodization. 
Concerning the transmission phase, the effect of the parameters is summarized by the transmit beam pattern, which encodes the information about energy distribution across the image field. 
While the large majority of ultrasound modalities deals with wideband beam patterns, the quantitative application Acoustic Radiation Force Impulse (ARFI) elastography employs narrowband beam patterns. 
ARFI is a specific elastography technique in which the probe emits a narrowband shock pulse to generate a shear wave in the tissue; several tracking pulses are then emitted to sample the propagation of the shear wave, and to estimate the local shear wave velocity; such an estimation allows one to compute the tissue elasticity map, displayed as a color-coded image~\cite{Deffieux2012,Kijanka2018,Hollender2015,Doherty2013,Rouze2012}.
Generally, the ideal transmit beam pattern, both narrowband or wideband, should have a main lobe width as uniform as possible along the depth of interest and a side lobe level as low as possible.  
Moreover, the main lobe intensity along depth should be as uniform as possible, since its maximum peak is bounded by both mechanical and thermal limits assuring the safety of the insonification.
These requirements are in contrast with the typical hourglass-shaped beam obtained with a fixed focalization depth~\cite{Jensen2002,Harris1981}, in which the main lobe width has its minimum nearby the focal depth and it widens before and above.  
Furthermore, the near field region below the focal depth is often characterized by side lobes whose intensity is even higher than the value along the central line. 
Finally, the combined effect of geometric and thermoviscous attenuation causes a quick energy drop-out immediately after the focal depth.  
These drawbacks have a negative impact on both image quality and measurement accuracy. 
For example, in ARFI elastography a non uniform main lobe width generates shear waves of non uniform temporal length along depth, causing bias in the estimation of shear wave times of arrival. 
The presence of side lobes with high intensity generates instead spurious shear waves that hinder the effect of the desired one, while the quick energy drop-out after the focal depth causes a rapid decrease of shear wave energy along depth, thus dramatically decreasing the SNR and the accuracy of the estimated shear wave time of arrival~\cite{Zhao2011}. 
Tuning the high level parameters of a standard beam pattern, for instance focal depth, aperture, and transmit frequency, may improve some beam pattern features but often at the expense of others. 
For example, increasing the aperture and moving the focal depth forward reduces the energy drop out but decreases the main lobe width uniformity.  
A different technique, named synthetic transmit beam~\cite{Mo2020} improves the transmit beam uniformity retrospectively, i.e., recombining different laterally spaced received signals.  
Though quite effective, this technique is not feasible for ARFI methods in which only one narrowband transmission is performed.
A different approach addressed in the literature~\cite{Kuperman1990,Curletto2003,Curletto2007,Trucco2002,He2015,Cardone2001,Jeong200} consists in optimizing a set of low level parameters, such as the transmit intensity at each channel, i.e., the apodization window, by minimizing a suitable cost function in the wideband case. 
In this framework the solution of this problem can have a closed form or may be easily computed with optimization algorithms, being the cost function convex.
A few other methods optimize the delay patterns solely regarding the intensity at focus~\cite{Herbert2009,Marsac2012,Marsac2010}, while other strategies aim at calculating the weight vector of a transducer array to obtain a low-sidelobe transmitting beam pattern~\cite{He2014,He2008}. To the best of our knowledge the optimization of transmit delays has been relatively underexplored, with a greater emphasis placed on advancing image reconstruction techniques. Meanwhile, the challenges associated with parameter tuning are presently addressed through manual tuning methods.

In this work we propose to optimize the narrowband beam pattern, suited for ARFI elastography, by considering each single transmit delay as an independent variable and keeping the transmit intensity fixed. 
The choice of optimizing apodization-free beam patterns, in which all the active elements transmit at the same voltage, widens the scope of application of this method also to low-price commercial apparati, equipped with simple two level transmit pulsers instead of linear amplifiers.
The much higher number of degrees of freedom, in comparison with the classical single-focus delay parametrization, allows the energy to be distributed over a greater depth range, thus improving the beam uniformity along depth. 
To compare quantitatively the optimized beam patterns with a large set of standard ones, we define three novel ad-hoc quality metrics, specifically designed to measure how much a beam pattern approaches the ideal behaviour required by the applications described above.
Existing beam pattern metrics, like peak side lobe level, integrated side lobe level, or full width half maximum of main lobe are quite effective for far field beam patterns or to evaluate near field beam pattern at a single depth. 
However, they fail to capture the global beam pattern quality along a depth range. 
Our numerical experiments demonstrate that the optimized BPs are indeed thinner and more uniform along depth, with low side lobes and a uniform main lobe intensity. The application of the method to ARFI elastography results in higher intensity and uniformity for shear waves, which will ultimately lead to more precise reconstruction of elasticity maps". It should be emphasized, however, that we hereby only consider the optimization of transmission step related to the shock pulse, leaving the tracking transmission and reception  for future studies.

The paper is organized as follows: Section~\ref{sec:materials} briefly describes the mathematical formulation of the problem including the derivation of the transmit BP equation, the definition of the optimization problem and of the metrics for evaluating the goodness of the results. 
Section~\ref{sec:results} shows the simulation setting and presents both the qualitative and the quantitative results.
Finally, in Section~\ref{sec:conclusions} our conclusions are drawn.
\section{Materials and Methods}\label{sec:materials}
\subsection{Biomedical ultrasound model formulation}
We proceed with describing the simulation of the transmit beam pattern (TBP) and we focus more specifically on the narrowband case, which allows to approximate the waveform as an infinite pulse. 
We refer to~\cite{Razzetta2023} for a full explanation of the rationale behind this reasoning, and we will give here only a brief description.
We consider a linear probe composed of $2N$ piezoelectric elements of width $\Delta x$ that convert the electric pulse into a pressure wave and vice-versa. For convex probes the treatment is analogous.
We consider a Cartesian coordinate system whose origin corresponds to the geometrical center of the probe (see Figure~\ref{fig:scheme}), so that the center of each probe element is at the position $c_i = (i + 1/2) \Delta x$ with $i=-N,\ldots,N-1$.
Each point $\vec{x} = (x,z)$ of the image plane is identified by a positive value of the $z$ coordinate.
\subsubsection{The BP equation}
The TBP is the pattern of radiation generated from a beam and it can be described as the energy or the power, for a periodic signal, of the propagating wavefront, evaluated at every spatial point of the field. 
The propagating wavefront can be computed from the signal emitted by the probe and the spatial-temporal impulse response characterising the probe itself. 
In particular, the impulse response function $h_i(\vec{x},t)$ of each $i$-th element is the integral over the element surface $A_i$ of a spherical wave starting from each point of the element surface~\cite{Harris1981}:
$$h_i(\vec{x},t) = \int_{A_i} \frac{\delta\left(t - \frac{\vert\vert \vec{x}-\vec{y}\vert\vert}{c}\right)}{2\pi\vert\vert \vec{x}-\vec{y}\vert\vert} d\vec{y},$$
where $c$ is the ultrasound speed propagation.
We assume that every active element emits the same temporal transmit waveform $I(t)$~\cite{Jensen1996}. 
Thus, the resulting signal is obtained at every point $\vec{x}$ by convolving $I(t)$ with the sum of all the impulse responses of the active elements and with a term accounting for the tissue thermoviscous attenuation
$a(\vec{x},t)$. 
The waveform $I(t)$ can be trasmitted with a different delay for each active element. 
The BP shape depends on the set of transmission times because of the various destructive and constructive interferences which are generated during the propagation.
We introduce the transmission delay $D_i$, associated with the $i$-th element, and we describe the delay in the signal emission as a translation in time $\delta(t - D_i)$. By defining 
$$H_i(\vec{x},\cdot) : = a(\vec{x},\cdot) \ast h_i(\vec{x},\cdot),$$
we can describe the propagating wavefront in time as:
$$S(\vec{x},\cdot) = I \ast \sum_{i=-N}^{N-1} H_{i}(\vec{x},\cdot) \ast \delta(\cdot - D_i).$$
In the narrowband case one can assume $I(t)$ to be an infinite pulse with given frequency $f_0$, i.e., $I(t) = e^{-2 \pi j f_0 t}$.
Thus, the power of the BP is
\begin{equation}
	P(\vec{x}) : =  \frac{1}{T}\int_{[0,T]} \left \vert S(\vec{x},t) \right \vert ^2\, d t \label{eq:energy},
\end{equation}
where $T = \frac{1}{f_0}$.
Moving to the time frequency domain and relying on Plancherel's theorem, stating the equivalence of signal power in time and frequency domain, the BP power can be written as the sum of the squared absolute values of the Fourier series coefficients:
\begin{equation}
	P(\vec{x}) = \sum_{s = -\infty}^{+\infty} \Big\vert \sum_{i=-N}^{N-1} \hat{H}_{i}(\vec{x}, s) e^{-2\pi j D_i f_0 s} \Big \vert^2, 
	\label{eq:coeffs}
\end{equation}
where $\hat{H}_{i}(\vec{x}, s)$ is the Fourier Series of $H_{i}(\vec{x}, t)$.
Since we have assumed $I(t)$ to be a pure tone, only one coefficient of the series is different from zero.
Therefore, Eq.~\ref{eq:coeffs} reduces to: 
$$P(\vec{x}) =  \Big \vert \sum_{i=-N}^{N-1} \hat{H}_{i}^{f_0}(\vec{x}) e^{-2\pi j D_i f_0} \Big \vert^2,$$
where we denote $\hat{H}_{i}^{f_0}(\vec{x}) : = \hat{H}_{i}(\vec{x}, 1)$.
In particular, we can exploit the modulus properties obtaining:
$$P(\vec{x}) = \left(\sum_{m=-N}^{N-1} \hat{H}_{m}^{f_0}(\vec{x}) e^{-2\pi j D_m f_0}\right) \overline{\left(\sum_{n=-N}^{N-1} \hat{H}_{n}^{f_0}(\vec{x}) e^{-2\pi j D_n f_0} \right)}.$$
Finally, making use of the exponential trigonometric form and using the cosine addition formula, the BP takes the form 
\begin{equation}
	P(\vec{x}) = \sum_{m,n=-N}^{N-1} {\hat{H}}_{m,n}^{f_0}(\vec{x}) \cos(2\pi f_0 (D_{m} - D_{n})) \label{eq:NB}
\end{equation}
where 
$$\hat{H}_{m,n}^{f_0}(\vec{x}) : = \hat{H}_m^{f_0}(\vec{x}) \cdot \overline{\hat{H}_n^{f_0}(\vec{x})}$$ 
is obtained as the product between the Fourier Series of $H_m(\vec{x}, t)$ and the conjugate of the Fourier transform of $H_n(\vec{x},t)$ at frequency $f_0$.

\subsection{Optimization of Beam Pattern}
This subsection introduces the optimization problem and three novel metrics that are based on the desired outcome for the application of the method to ARFI elastography.
\subsubsection{Delay optimization}
\begin{figure}
	\centering
	\includegraphics[width=\textwidth]{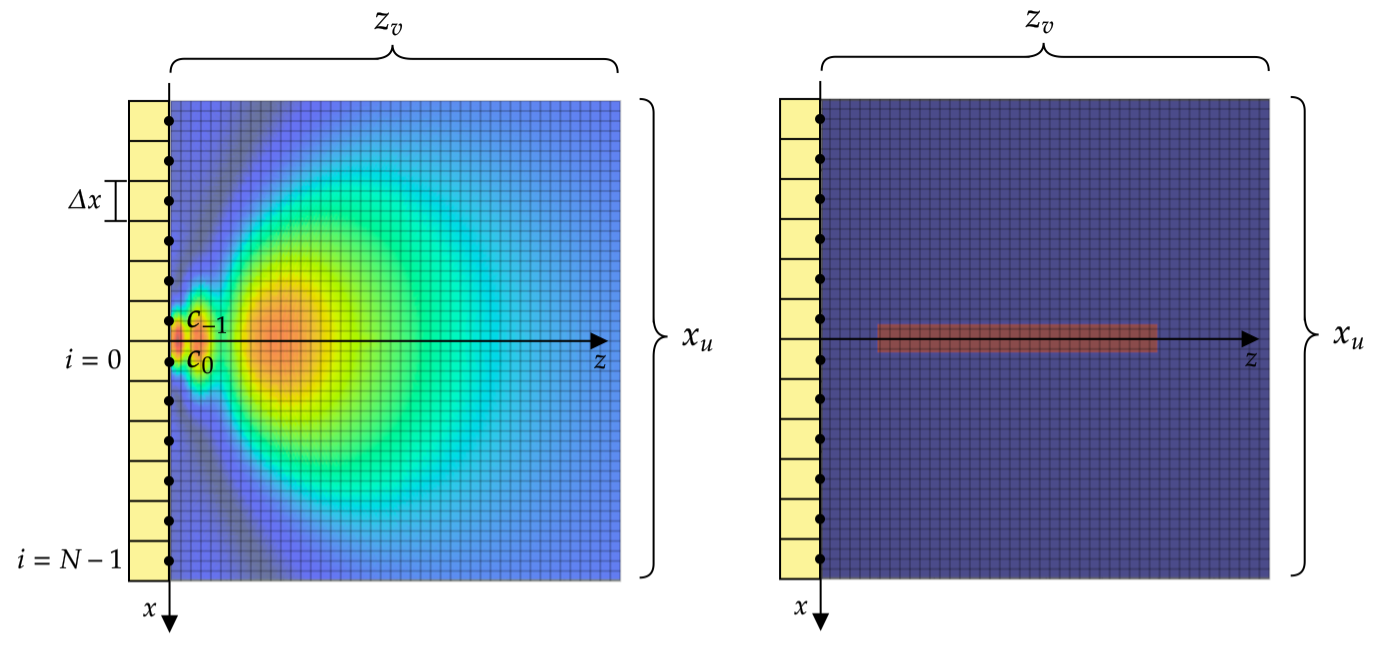}
	\caption{Left: Scheme of the adopted geometry, in the background the simulated energy pattern for elements $c_0$ and $c_{-1}$. Right: Prescribed rectangular BP to approximate.} \label{fig:scheme}
\end{figure}

The aim of the optimization problem is to make $P(\vec{x})$ as similar as possible to a prescribed BP shape by varying the frequency $f_0$, the number of active elements $N$, and their corresponding delays $D_i$, with $i=0,\dots, N-1$.
Since we are interested in BP shapes symmetric with respect to the probe central transmit line, we impose $D_{-i}=D_{i-1}$. 
We consider the delays as free variables so that the BP changes according to $N + 2$ parameters:
$$
p := (N, f_0, D_0, \dots, D_{N-1}).
$$
In order to compute the BP in Eq.~\ref{eq:NB}, we sample $\vec{x}$ by taking a grid of points $(x_u,z_v)$ with $x_u = u \Delta x/4$ where $u=-4N,\ldots,4N$ and $v=0,\ldots,M$ (see Figure~\ref{fig:scheme}). 
From a computational point of view, we approximate $h_i(x_u,z_v)$, the integral over a single element surface, according to~\cite{Jensen2002,Jensen1996}.
Moreover, by fixing the central frequency $f_0$, we pre-compute the attenuation coefficients and the impulse response function for each element. 
An example of impulse response function for a linear probe with pulse central frequency $f_0=4.5$ MHz is given on the left panel of Figure~\ref{fig:scheme}.
Then we compute $\hat H_{i, j}^{f_0}(x_u,z_v)$ by using the discrete Fourier transform over time at frequency $f_0$.
The sampled BP is a vector $X(p)$ depending non-linearly on $p$, whose components are given by $X_k := P(\vec{x}_k)$ where $k$ is an unrolled index, namely $k=u + (8N+1) v$.
Therefore, we can formulate a Least Squares (LS) problem as follows
\begin{equation}
	p^* = \arg \min_{p} \Vert X(p) - G \Vert^2_2, \label{eq:LS}
\end{equation}
where $G$ is the discretization of the prescribed BP rectangular shape we want to approximate. 
The width of $G$ has been chosen (experimentally) to be $4$ pixels ($=\Delta x$), i.e., the minimum width possible so that the optimization algorithm would give non trivial results (limit case of no transmission).
Since the vector $X(p)$ is separately periodic with respect to the delays $D_0, \ldots, D_{N-1}$, the optimization problem is non linear and non convex as the parameter domain contains an $N$-dimensional torus (cf. Appendix in~\cite{Razzetta2023}), therefore not suitable for standard optimization techniques.
Among the many methods developed to deal with non convex optimization, e.g. genetic algorithms, simulated annealing, manifold optimization to mention a few, we propose to use the Particle Swarm Optimization.
This method has been proven to be particularly efficient in solving this problem, after a number of attempts with different methods.

\subsubsection{Evaluation Metrics}
\label{ssec:EvaluationMetrics}
\begin{figure}[!ht]
	\centering
	\includegraphics[width=\textwidth]{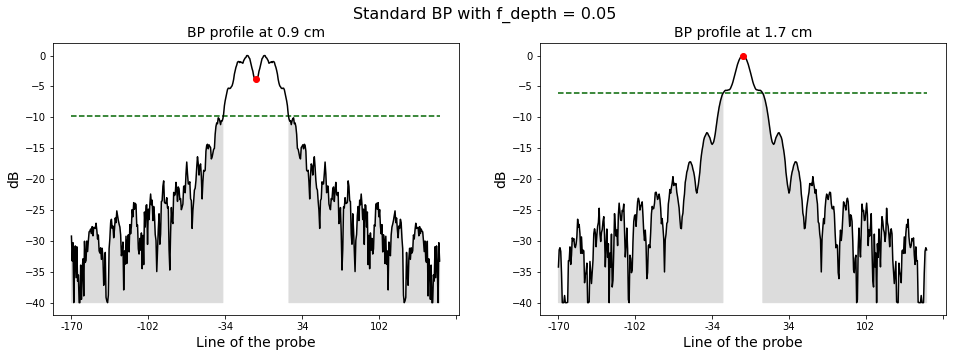}
	\caption{Two examples of BP profile to clarify the rationale for the normalization with respect to the central line in the metrics definition. The red dot identifies the BP central line, the green dashed line highlights the -6 dB level under the central value, the grey area is the Side Lobe area and the central white area is the main lobe area. The main lobe width increases when the maxima are not on the central line, thus penalizing BPs with this undesired behaviour.}
	\label{fig:immmis}
\end{figure}

In order to allow a quantitative evaluation of optimized BPs with respect to standard ones, we devised a set of metrics aimed at providing a synthetic evaluation of the overall BP quality. 
The first metric aims at measuring the main lobe width average value and uniformity along depth.
To do so, we make the assumption that the Main Lobe (ML) at each depth consists of the BP portion centered at the transmission line that differs from the central value less than $6$ dB (Figure~\ref{fig:immmis}). 
The choice of this threshold is coherent with the usual definition of Full Width at Half Maximum. Moreover, this assumption penalizes those patterns whose maximum value is not on the central line.
Therefore, we define the Main Lobe Width ($MLW$) at a certain fixed depth $z$ as:
\[MLW(z) = \max_{x,y \in \Psi_z}\vert x-y\vert \quad \textrm{ where } \quad \Psi_z := \left\{ x \quad \text{s.t.} \quad P(x, z) \geq \frac{P(\bar{x}, z)}{\epsilon} \right\},\]
where $\bar{x}$ is the azimuth coordinate of the beam central line and $\epsilon = 10^{\frac{3}{5}} \approx 4$ .

The second metric evaluates the overall side lobe level. 
Let us define the Side Lobe Level ($SLL$) at a fixed depth as the average over $x$ of all the BP values outside the region of the main lobe (as defined above): 
$$SLL(z) = \frac{1}{P(\bar{x},t)} \frac{1}{|\Psi_z^C|} \int_{\Psi_z^C} P(x,z) \, dx ,$$
where $|\Psi_z^C|$ denotes the measure of the complementary set $\Psi_z^C$ of  $\Psi_z$.
Note that the $SLL$ is normalized with respect to the central line value at the each depth $z$, in order to cope with the overall decrease trend of the BP along depth.
By computing mean and standard deviation of $MLW(z)$ and $SLL(z)$ over a set of depths of interest, we can evaluate the main lobe width uniformity and narrowness, and the side lobes uniformity and average level.

Finally, to evaluate the drop out of main lobe power along depth we compute the mean integral of the central line values, namely the Central Line Power ($CLP$):
$$CLP = \frac{1}{|\zeta|}\int_{\zeta} \frac{P(\bar{x}, z)}{\max_{\hat{x}, \hat{z}}P(\hat {x}, \hat{z})}   \, dz \quad \textrm{ where } \quad \zeta = \left\{z \quad \text{s.t.} \quad z \in [z_{\min}, z_{\max}] \right\},$$
where $[z_{min}, z_{max}]$ denotes a range of depths of interest.
This quantity measures how uniform the line is along depth and it gives information about the power drop.

\section{Numerical results and discussion}\label{sec:results}
This section presents our numerical experiments. 
We start by briefly explaining how realistic measurement data are simulated and by reviewing the coding strategy that is based on a greedy approach and Particle Swarm Optimization (PSO); for more information consult, e.g.,~\cite{Kennedy1995}. Finally, some qualitative and quantitative results are presented.
\subsection{Simulation setting and coding strategy}\label{ssec:setting}
Let us consider two linear probes (Probe $1$ and Probe $2$) with the same depth range of $4.5$ cm and different pitch $\Delta x$ (smaller for Probe $2$). 
Our purpose is to obtain a beam shape as uniform as possible along depth, with the lowest possible level of side lobes. 
The goal is to start from the usual hour-glass shaped beam pattern related to a fixed focal depth, by focalizing each element at a different depth to achieve a more stretched shape.
We define the desired BP as a rectangular shape centered at depth $1.75$ cm with length of $2.5$ cm, with value one inside and zero outside the rectangle (right panel of Figure~\ref{fig:scheme}). 
The choice of this particular rectangular shape is reasonable under different aspects. 
First, it forces the BP uniformity along depth, while the zero values at its sides enforce a low level of side lobes. 
Also, it should be noted that the rectangle does not fully cover the whole depth range, so as to disregard the BP behavior before the first depth of interest, where energy levels are consistently high, as well as beyond a specific depth. 
This is due to the attenuation factor hindering signal transmission with high frequencies, which makes optimizing the final field portion challenging for linear probes.

As for the optimization algorithms, we employ a greedy approach to optimize both the frequency of the narrowband pulse and the aperture, while we exploit Particle Swarm Optimization (PSO) Algorithm~\cite{Kennedy1995} for optimizing the vector of delays, with the following rationale. By choosing a finite set of frequencies suitable for the probe, we can generate the corresponding impulse response map set and rapidly simulate a beam. 
The aperture is an integer and it is reasonable to assume the active surface length is greater than the smallest dimension of the rectangle. 
The choice of PSO algorithm is justified by the nature of the problem: the free parameters' domain contains a N-dimensional torus (cf.~\cite{Razzetta2023} for the reasoning behind this assertion) so every gradient-based algorithm would work in a fixed local chart, thus being very sensitive to the initialization. 
Conversely, PSO adopts multiple random initialization hence being less affected by the manifold domain structure.

To generate the beam images we use our simulator parUST~\cite{parust}, developed in Python\texttrademark and freely available on GitHub. 
The simulator implementation is based on a two-step technique that allows one to pre-compute the system frequency response for a given probe model and region of interest: this methodology is particularly suited when the simulator is used as a kernel in gradient-based or non gradient-based optimization problems to the effect that many different sets of transmit parameters on the beam pattern shape have to be tested.
Finally, we usePySwarms\footnote{https://pyswarms.readthedocs.io/en/latest/} free software research toolkit to apply PSO algorithm in Python\texttrademark. 
At each iteration the algorithm generates a BP for each particle corresponding to a vector of transmit delays, and it measures the adherence between the actual and the prescribed beam. Afterwards, it updates each particle by taking into account its adherence value and the overall direction of the swarm. 
Regarding the algorithm parameters, we set the cognitive parameter to $0.5$, the social parameter to $0.7$, and the inertia parameter to $0.9$, thus giving more importance to the group behaviour rather than to the single particle movement.
It should be pointed out that the overall optimization computational cost can be high (around $8$ hours for one experiment), but, on the other hand, the idea is to provide an initial suboptimal solution for conducting practical tests with the machines and then build upon them using a solid set of parameters.

\subsection{Qualitative results}
\begin{figure}
	\centering
	\includegraphics[width=\textwidth]{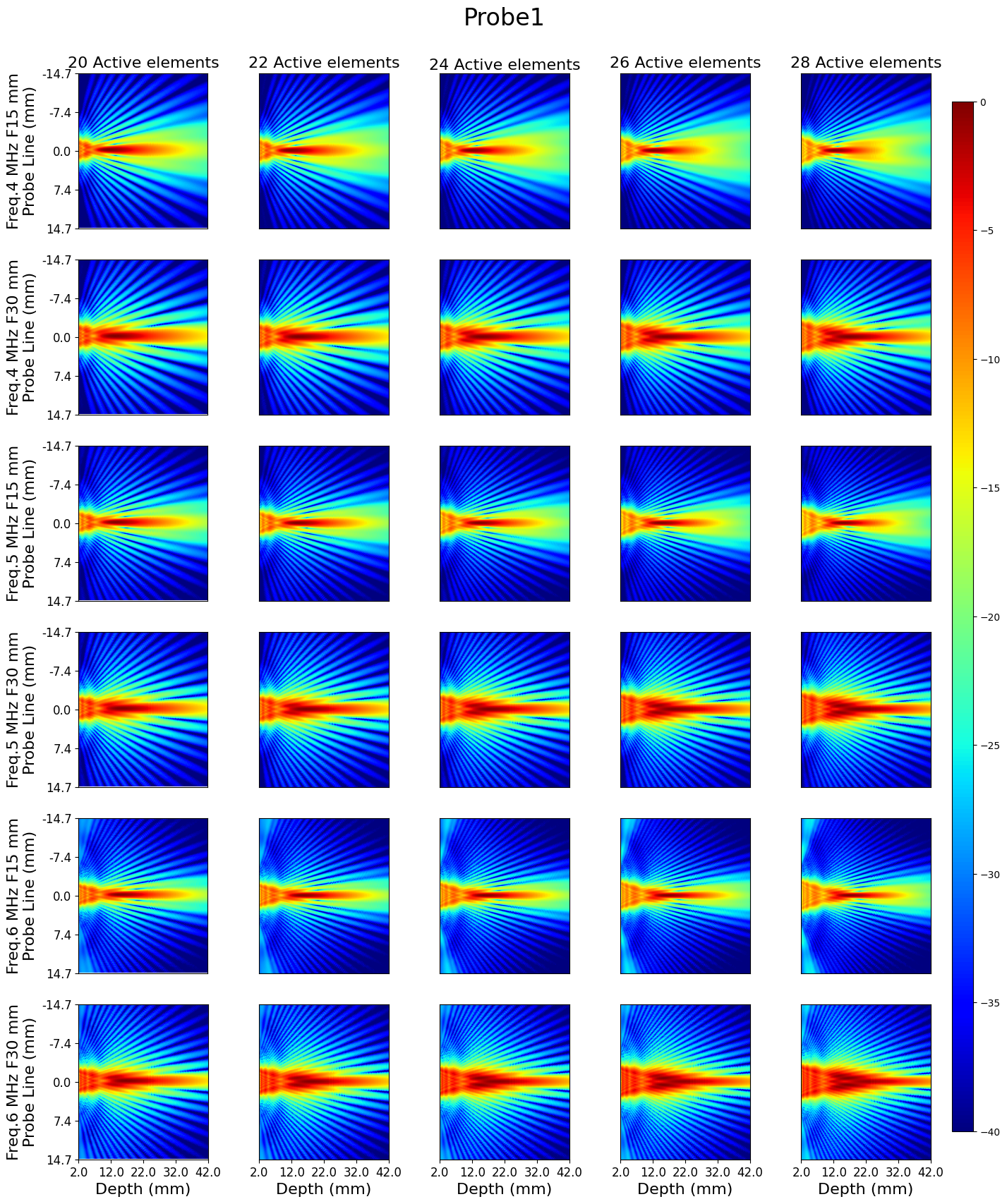}
	\caption{Narrowband beam shapes obtained with standard focal law (single depth focus) at frequency $4.5$ MHz with Probe 1. The focal depth is constant along the rows, while the aperture is constant along the columns.} \label{fig:std1}
\end{figure}

\begin{figure}
	\centering
	\includegraphics[width=\textwidth]{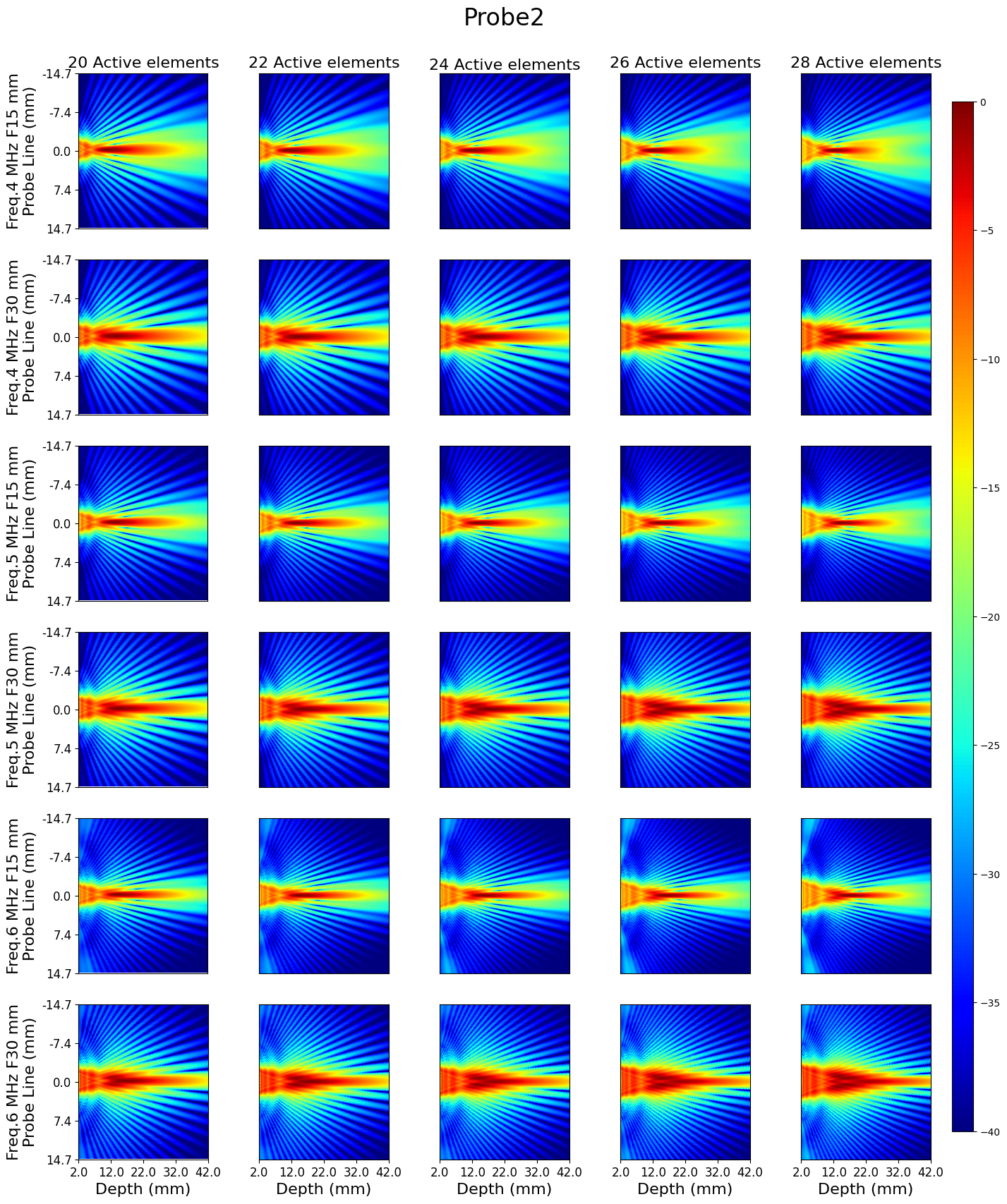}
	\caption{Narrowband beam shapes obtained with standard focal law (single depth focus) at frequency $4.5$ MHz with Probe 2. The focal depth is constant along the rows, while the aperture is constant along the columns.} \label{fig:std2}
\end{figure}

\begin{figure}
	\centering
	\includegraphics[width=\textwidth]{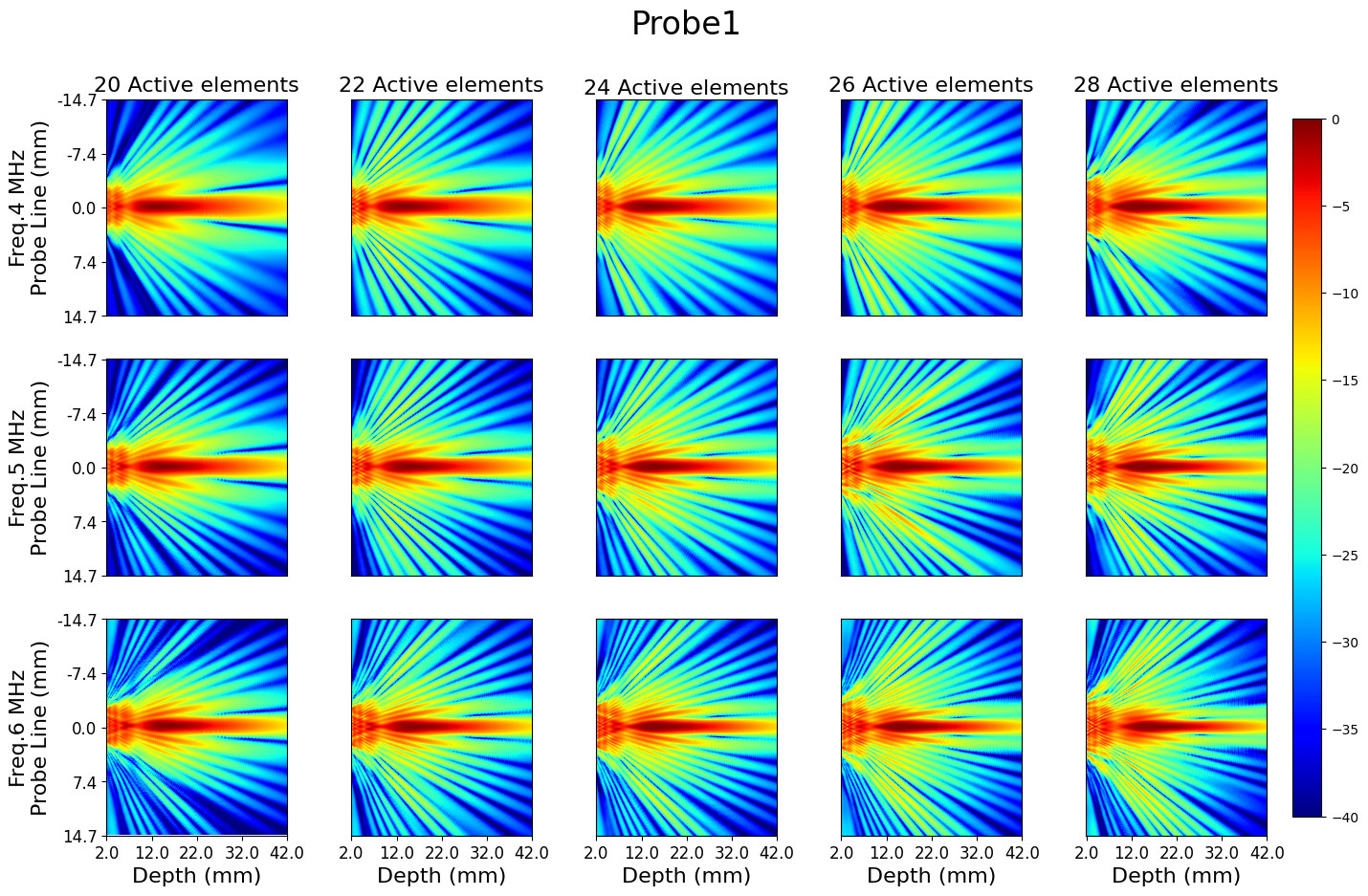}
	\caption{Optimized beam patterns for Probe 1 for three different frequencies and with five different numbers of active elements.}
	\label{fig:resultsPB1}
\end{figure}

\begin{figure}
	\centering
	\includegraphics[width=\textwidth]{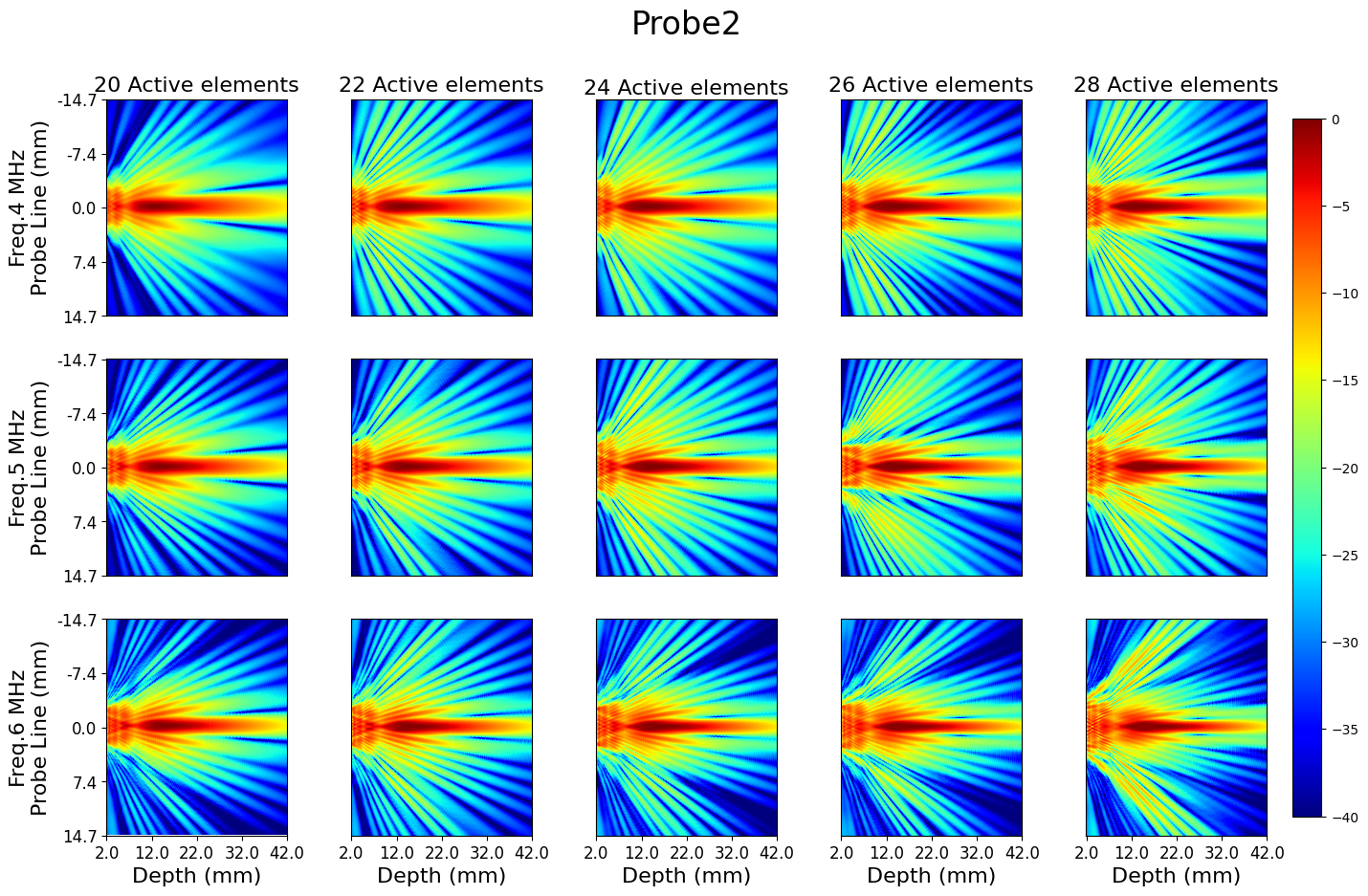}
	\caption{Optimized beam patterns for Probe 2 for three different frequencies and with five different numbers of active elements.}
	\label{fig:resultsPB2}
\end{figure}

First, we qualitatively analyze a set of narrowband BPs obtained with a standard focalization law as a function of focal depth, aperture, and frequency (as displayed in Figures~\ref{fig:std1} and~\ref{fig:std2}), and compare their shapes with the optimized beam patterns for the two different probes described in~\ref{ssec:setting} (Figures~\ref{fig:resultsPB1} and~\ref{fig:resultsPB2}). 

The well-known drawbacks of standard delay profile can be observed in Figures~\ref{fig:std1} and~\ref{fig:std2} for Probe 1. 
For small apertures, the main lobe is considerably large and its intensity profile rapidly decreases over depth. 
By increasing the aperture, the main lobe width shrinks around the focal depth but it is enlarged considerably in the near field region, where the intensity peak is not even aligned with the main lobe axis.
These drawbacks can severely affect the accuracy of shear wave propagation estimation in ARFI elastography~\cite{Zhao2011}. 
The same holds for Probe 2, with the only difference being that the probability of producing grating lobes is reduced, primarily because of the smaller pitch.

Figures~\ref{fig:resultsPB1} and~\ref{fig:resultsPB2} display some examples of the BP obtained by optimizing aperture, frequency, and delay profile.  
Focusing the analysis on the depth range of interest given by the rectangle of the prescribed beam pattern (right panel in Figure~\ref{fig:scheme}), it can be observed that the intensity peak is located on the main lobe axis at each depth. 
Moreover, the main lobe width is considerably uniform over the whole depth range especially in the near field in comparison with the BP obtained with the same aperture ($26$ elements) and a standard delay profile. 

\begin{figure}[!ht]
	\centering
	\includegraphics[width=\textwidth]{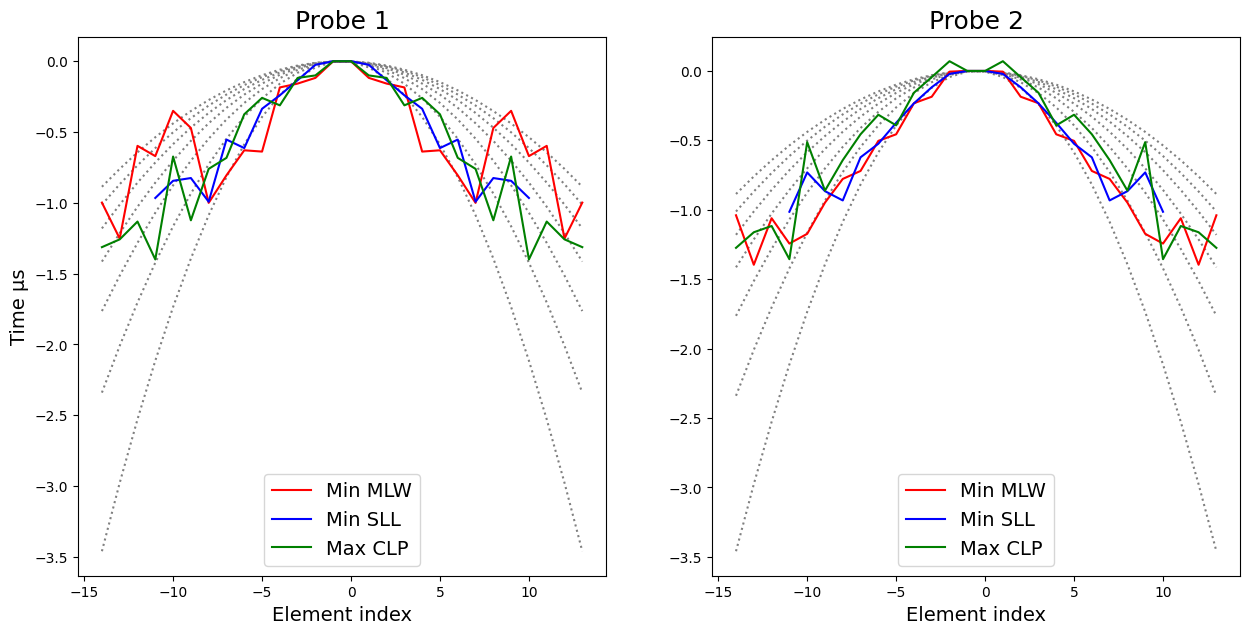}
	\caption{Comparison between standard and optimized delay profiles for both Probe 1 and 2. Standard profiles are depicted with the dashed grey curves and they get higher with increasing focal depth (from $0.5$ to $3.5$ $cm$). Optimized curves are depicted in red, blue, and green, corresponding to minimum $MLW$ value, minimum $SLL$ value, and maximum $CLP$ value, respectively (see Section~\ref{ssec:EvaluationMetrics}).} 
	\label{fig:delays}
\end{figure}

These improvements are achieved thanks to the optimized delay profile, displayed in Figure~\ref{fig:delays} together with the standard profiles for three different focal depths. 
The relationship of the optimized profile to the family of standard profiles confirms the effectiveness of the proposed approach in optimizing the narrowband TBP.
Considering each delay as a free variable has increased significantly the number of delays combinations, therefore augmenting the possible BP shapes. 
In particular, the optimized profile can be seen as a combination of standard delay profiles at different focal depths.

\subsection{Quantitative results}

\subsubsection{Evaluation Metrics outcomes}

\begin{figure}[!ht]
	\centering
	\includegraphics[width=\textwidth]{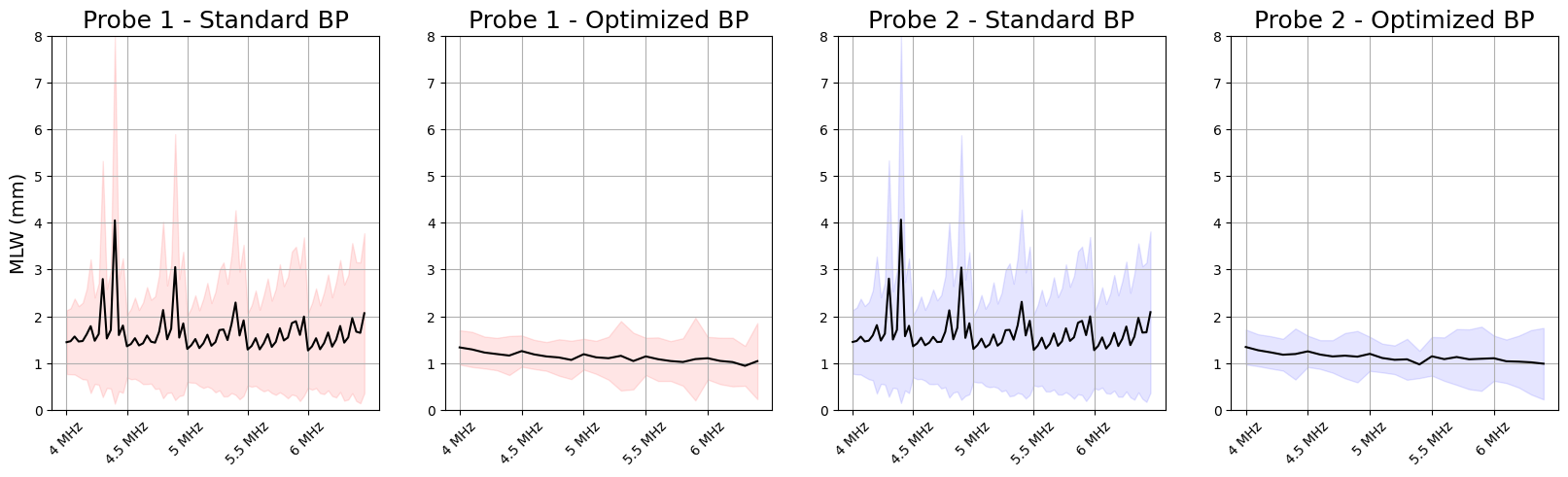}
	\caption{Values of $MLW$ (mm) over 2 probes (red and blue), 5 frequencies, 5 apertures. For standard BPs, $3$ focus values are considered ($1.5, 2, 3$ cm), hence $3$ points in the standard BP plot corresponds to 1 point in the optimized version.} \label{fig:MLW}
\end{figure}

\begin{figure}[!ht]
	\centering
	\includegraphics[width=\textwidth]{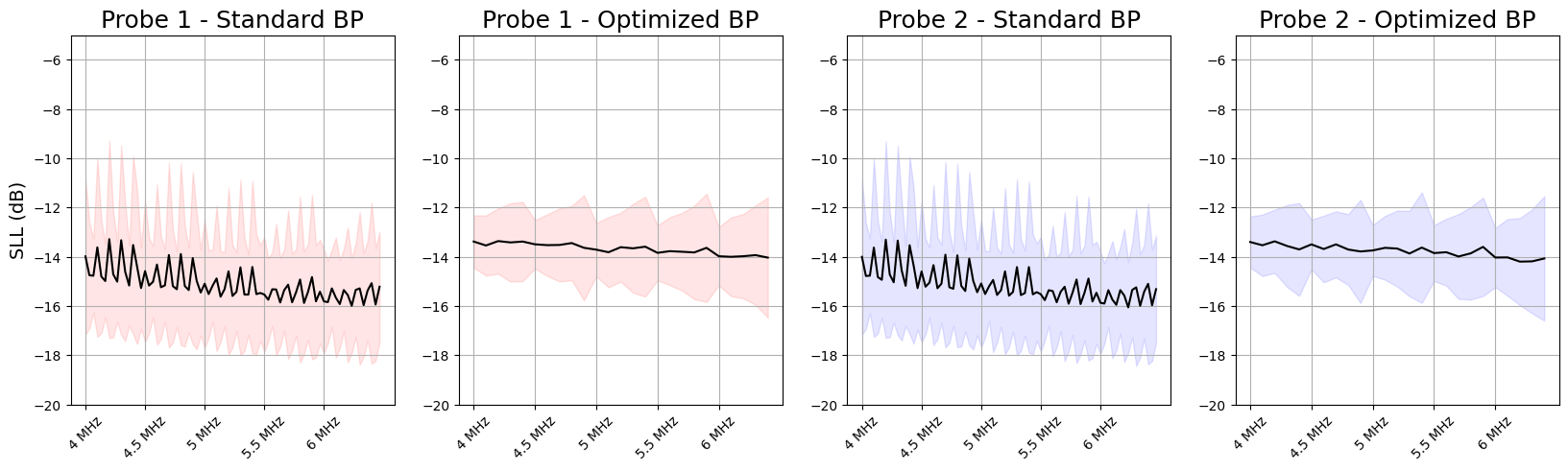}
	\caption{Values for $SLL$ (dB) over 2 probes, 5 frequencies, 5 apertures.  For standard BPs, $3$ focus values are considered ($1.5, 2, 3$ cm), hence $3$ points in the standard BP plot corresponds to 1 point in the optimized version.} \label{fig:SLL}
\end{figure}

\begin{figure}[!ht]
	\centering
	\includegraphics[width=\textwidth]{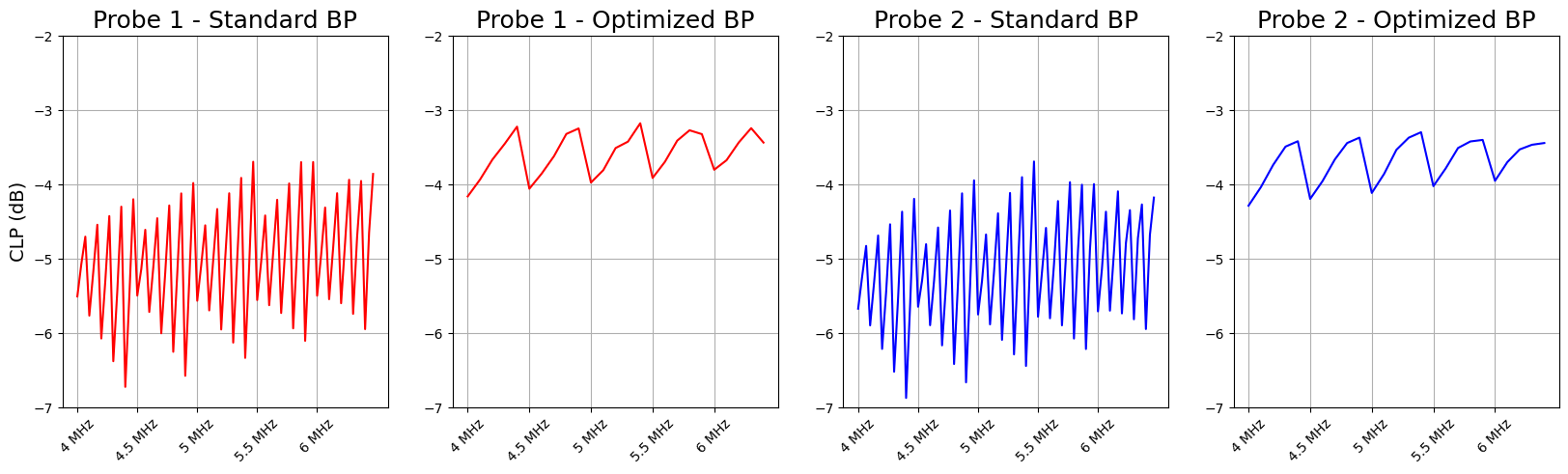}
	\caption{Values of $CLP$ (dB) over 2 probes, 5 frequencies, 5 apertures.  For standard BPs, $3$ focus values are considered ($1.5, 2, 3$ cm), hence $3$ points in the standard BP plot corresponds to 1 point in the optimized version.} \label{fig:CLP}
\end{figure}

Figures~\ref{fig:MLW}, \ref{fig:SLL} and \ref{fig:CLP} provide a graphical depiction of values for $MLW$, $SLL$, and $CLP$, both in the standard case and for the optimized beam patterns, computed in the case of the above introduced probe types. 
In both cases, metrics values refer to numerical experiments with $5$ different frequency values and $5$ apertures, while for standard BPs also $3$ focus values are taken into account. 
Metrics values are computed along the depth interval defined from $z_{min} = 1.75$ cm until the end of the field of interest ($4.2$ cm). 
This indicates that the chosen rectangle lenght does not severely affect the outcome of the optimization algorithm, which remains better than the standard cases even outside the rectangle limits. 
It can be observed that $MLW$ values for the optimized BPs are always lower than the standard case, and with a significantly smaller standard deviation, in line with what one should expect from the optimization algorithm.
This indicates that all the $25$ optimized BPs have a thinner and more uniform $MLW$ with respect to all the $5\cdot 5\cdot 3 = 75$ evaluated standard BPs (cf. Figure~\ref{fig:MLW}).

As for $SLL$ values, they are on average slightly higher than the standard ones, but considerably more stable throughout all the experiments (cf. Figure~\ref{fig:SLL}).
Indeed, one can observe that $SLL$ for the optimized BPs achieve values in between the best and the worst values of standard BP, both for mean and standard deviation. 
However, the lowest average $SLL$ for standard BP is just 1.5 dB lower than the lowest average $SLL$ for optimized BP and is achieved at the price of a much higher $MLW$.

Finally, our investigation confirms that $CLP$ values are higher in the optimized case for both probes, once again in line with the desired criteria for optimization (cf. Figure~\ref{fig:CLP}).
More significantly, the best values of $CLP$ for standard BP are in general achieved with large apertures, but this condition is associated with worse values of $MLW$.

These results support the thesis that the proposed optimization method enforces the prescribed beam pattern behaviour also outside the optimization zone, thus making the choice of prescribed BP not critical. 

\subsubsection{ARFI evaluation}
In the context of ARFI elastography, a single ultrasound transducer serves the dual purpose of both inducing and monitoring a deformation response \cite{Nightingale2011}. 
The distribution of the acoustic radiation force field is spatially governed by the active transducer aperture, the material properties, and characteristics of the transmitted beam, while its magnitude is influenced by factors such as attenuation and intensity. 
Since attenuation varies with frequency and depth, selecting the optimal frequency to generate an acoustic radiation force becomes an application-specific task and entails a trade-off between minimizing attenuation losses in the near field and maximizing focal point gain ~\cite{Palmeri2006}. 
The configuration of the transmitted excitation beam parameters can be customized to account for these effects.

\begin{figure}
	\centering
	\includegraphics[width=\textwidth]{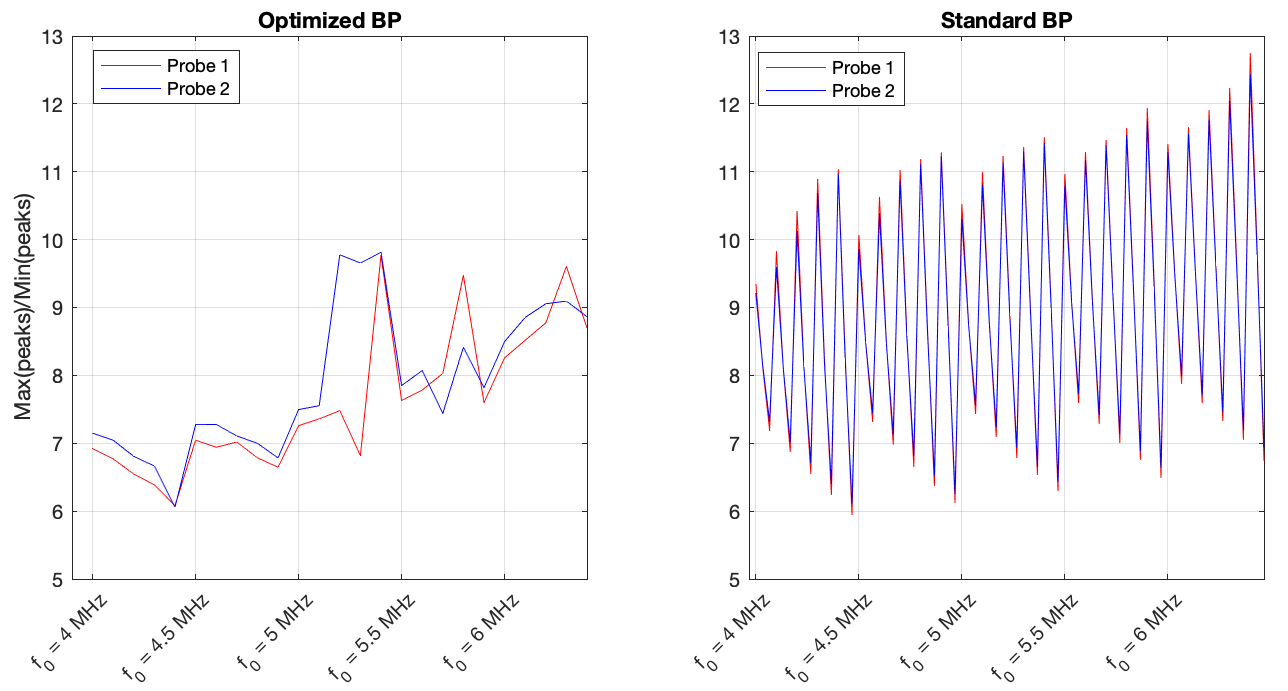}
	\caption{Range of shear wave peak intensities, the lower the range the better.} 
	\label{fig:MinMax}
\end{figure}

\begin{figure}
	\centering
	\includegraphics[width=\textwidth]{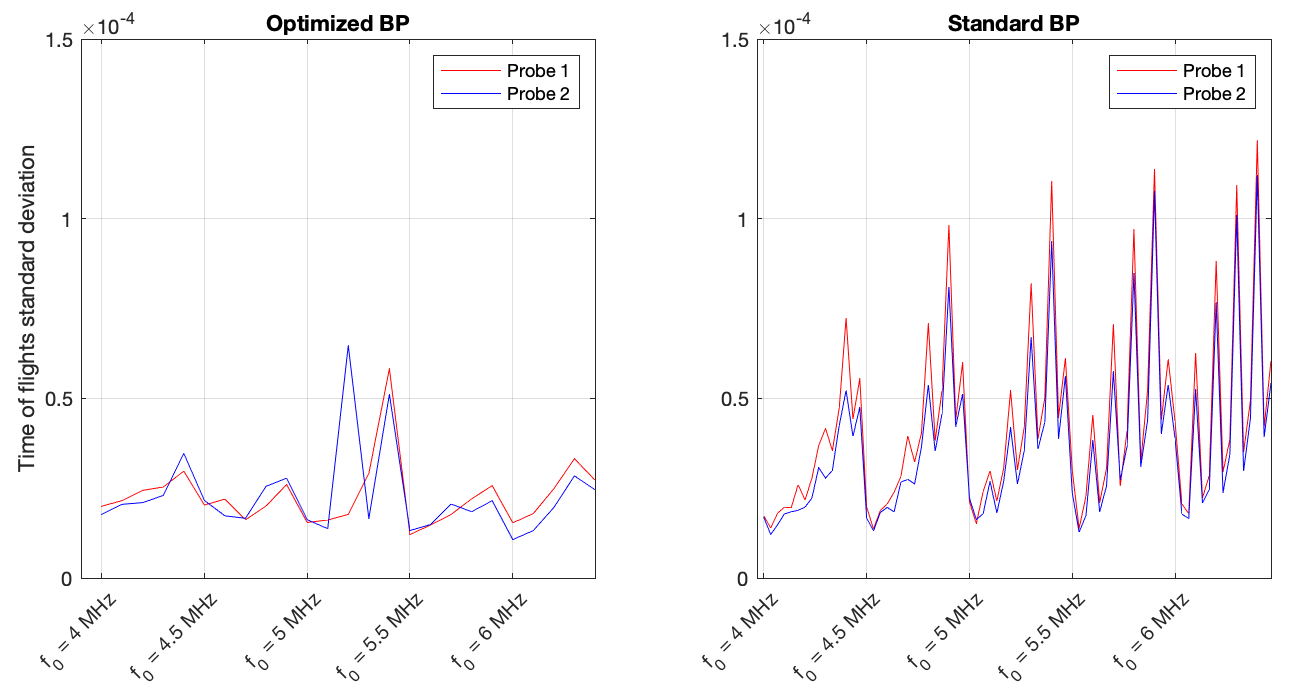}
	\caption{For each considered tracking distance from the center, we compute the standard deviation of times of flights along depth. We display the mean of these standard deviations. The lower the mean, the more uniform the time of flights.} 
	\label{fig:tflights}
\end{figure}

\begin{figure}
	\centering
	\includegraphics[width=\textwidth]{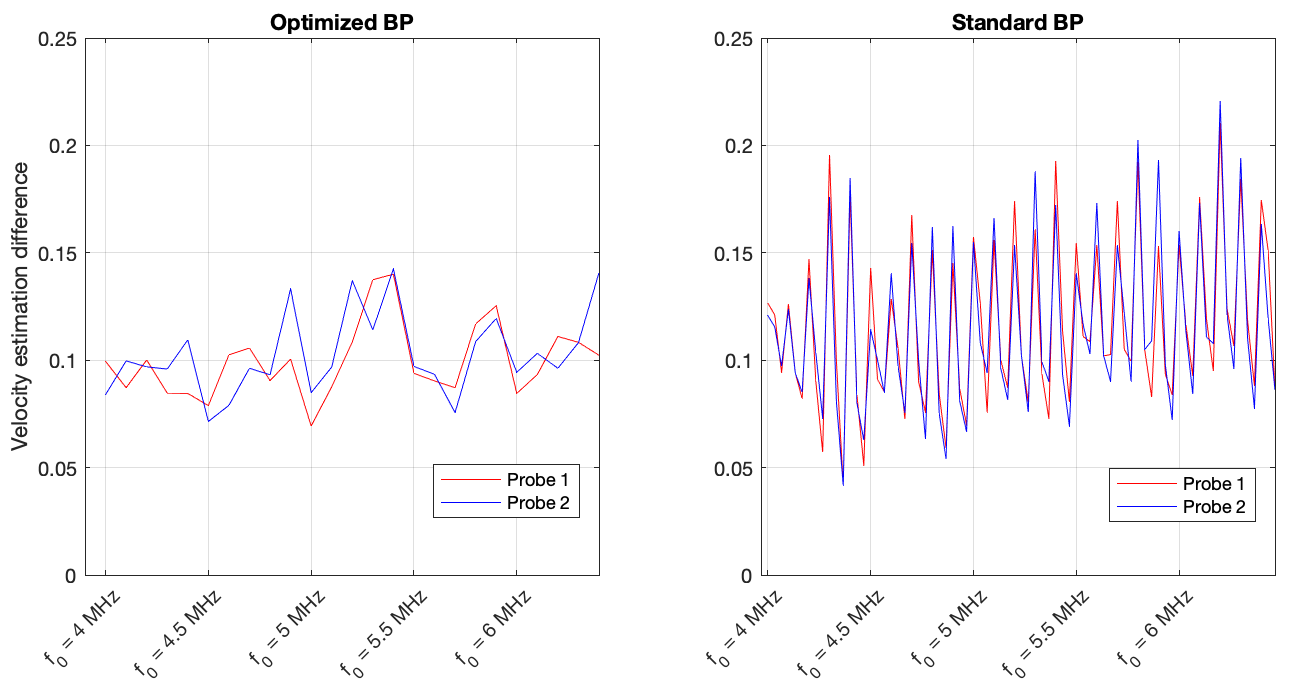}
	\caption{For each considered depth, we estimate the shear wave velocity along the tracking distance. In the graphs we report the absolute difference between the mean of the seven estimated velocities and the ``theorical" one. The lower this difference the better.} 
	\label{fig:vel}
\end{figure}

In order to assess that the optimization of the transmitted beam profile directly affects the desired characteristics of the generated shear waves, the following procedure has been implemented. 
First, the intensity of the transmitted field was simulated with a code based on Field II both with standard parameters and with optimized parameters.
Then, the intensity distribution of the transmitted fields was used as a shock source for shear wave generation~\cite{Huang2013}.
For each shock, we consider the same set of tracking coordinates where to estimate the shear wave.
We set as tracking coordinates an azimuth distance of $0.7, 1, 1.3, 1.5$ $cm$ from the shock pulse central line, and depth from the probe surface $0.5, 1, 1.5, 2, 2.5, 3, 3.5$ $cm$.
We evaluate the precision of ARFI estimation using three distinct measures. 
For every combination of depth and azimuth distance, we calculate 1) the shear wave time of flight, 2) maximum peak intensity, and for each depth we estimate 3) the shear velocity along azimuth performing a linear regression of times of flight.
A high peak intensity produces a higher Signal to Noise Ratio (SNR), which states a lower bound for the shear wave detection.
Because of this, the best condition is achieved when the  intensity of the peaks is uniform over all the tracking points.
For this reason we use as a quality metric the ratio between maximum and minimum values of the peaks.
The obtained results, shown in Figure~\ref{fig:MinMax}, demonstrate that optimized and standard BP are comparable under this criterion.
On the other hand, the uniformity of times of flight across depth fits the propagation model of an ideal cylindrical wave, which is at the base of shear velocity estimation methods. 
In Figure~\ref{fig:tflights} we show the mean standard deviation of times of flight in the optimized and the standard case, demonstrating the improvements achieved with the optimized BPs.
Lastly, we assess for every depth the precision of shear velocity estimation compared to the theoretical value imposed for simulation, and we display in Figure~\ref{fig:vel} the mean value along depth.
The obtained results confirm the optimized cases are at least comparable with the standard ones.


\section{Conclusions}\label{sec:conclusions}
Our work provides an effective approach for optimizing the narrowband transmit beam pattern, therefore providing a method that aims to ultimately improve accuracy and reduce artifacts in important medical ultrasound applications with a focus on ARFI elastography imaging.
It should be noted that we do not hereby consider the imaging improvements, but we provide all the necessary tools to quantitatively assess the strengths of our optimization method.
Indeed, the method provides sets of BPs that are extremely uniform along depth in terms of both BP width and intensity, and with low side lobe levels, as required by usual biomedical ultrasound applications. Moreover, the direct application to ARFI elastography provides comparable results in terms of shear wave peak intensity and even better results in terms of standard deviation of time of flights and velocity estimation.
The adopted strategy makes it easy to extend the calculations to convex probes,  by a simple change from Cartesian to polar coordinates, as well as to steered BP. 

Our work has some limitations nonetheless.
The advantages of this method cannot be proven by the optimized BPs with different number of active elements, nor with a different BP shape. 
In this regard, a potential expansion of the approach could involve experimenting Neural Networks with various ground truth shapes as inputs and the optimized delays as outputs. Future developments will also include apodization weights in the set of to-be-optimized variables, thus compensating the slight deterioration of side lobe level observed with the present method with respect to standard delay laws~\cite{Szabo2004}.
Moreover, the extension to wideband beam patterns will be investigated, as well as the extension of the method to imaging experiments.

\section*{Acknowledgment}
The authors wish to acknowledge the support obtained by Esaote S.p.A. under the research project ``MyLab 4.0", Code CUP\_B46G20001250005 funded by REACT EU (Operational Programme Enterprise and Innovation for Competitiveness 2014-2020). This work was carried out within the framework of the project “RAISE - Robotics and AI for Socioeconomic Empowerment” and has been partially supported by the European Union - NextGenerationEU. The views and opinions expressed herein are those of the authors alone and do not necessarily reflect those of the European Union or the European Commission. Neither the European Union nor the European
Commission can be held responsible for them

\printbibliography[heading=bibintoc] 
\end{document}